\documentclass{article}
\usepackage{amssymb,amsmath}
\usepackage{anysize}
\usepackage{psfrag}
\newtheorem{thm}{Theorem}[section]
\newtheorem{def.}{Definition}[section]

\newtheorem{prop}{Proposition}[section]
\newtheorem{cor}{Corollary}[section]

\newtheorem{lem}{Lemma}[section]

\numberwithin{table}{section}

\begin{document}
\title{The Teneva Game}
\author{Louis H. Kauffman\\
        Department of Mathematics, Statistics and Computer Science\\
        University of Illinois at Chicago\\
        851 S. Morgan St., Chicago IL 60607-7045\\
        USA\\
        \texttt{kauffman@uic.edu}\\
        and\\
        Pedro Lopes\\
        Center for Mathematical Analysis, Geometry and Dynamical Systems\\
        Department of Mathematics\\
        Instituto Superior T\'ecnico\\
        Technical University of Lisbon\\
        Av. Rovisco Pais\\
        1049-001 Lisbon\\
        Portugal\\
        \texttt{pelopes@math.ist.utl.pt}\\
}
\date{August 9, 2012}
\maketitle

\bigbreak

\begin{abstract}
For each prime $p > 7$ we obtain the expression for an upper bound on the minimum number of colors needed to non-trivially color $T(2, p)$, the torus knot of type $(2, p)$, modulo $p$. This expression is $t + 2 l -1$ where $t$ and $l$ are extracted from the prime $p$. It is obtained from iterating the so-called Teneva transformations which we introduced in a previous article.

With the aid of our estimate we show that the ratio ``number of colors needed vs. number of colors available'' tends to decrease with increasing modulus $p$. For instance as of prime $331$, the number of colors needed is already one tenth of the number of colors available.

Furthermore, we prove that $5$ is the minimum number of colors needed to non-trivially color $T(2, 11)$ modulo $11$.

Finally, as a preview of our future work, we prove that $5$ is the minimum number of colors modulo $11$ for two rational knots with determinant $11$.
\end{abstract}

\bigbreak


\bigbreak


\bigbreak

\section{Introduction.} \label{sect:intro}

\begin{figure}[!ht]
	\psfrag{a}{\huge$a$}
	\psfrag{b}{\huge$b$}
	\psfrag{c}{\huge$c=2b-a$}
	\psfrag{5}{\huge$\mathbf{5}$}
	\centerline{\scalebox{.5}{\includegraphics{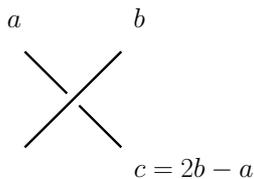}}}
	\caption{Colors of a coloring at a crossing.}\label{fig:x}
\end{figure}

\noindent
The reader will recall that a Fox coloring of a knot diagram is a labeling of the arcs of the diagram with the elements of $\mathbf{Z}/m\mathbf{Z}$ for some integer $m>1$, such that at each crossing the sums of the labels of the under-crossing arcs is equal to twice the labeling of the over-crossing arc. This is often expressed by the equation $c=2b-a$ where $a$ and $c$ are the labels for the under-crossing arcs and $b$ is the label for the over-crossing arc (see Figure \ref{fig:x}).

\begin{figure}[!ht]
	\psfrag{0}{\huge$0$}
	\psfrag{1}{\huge$1$}
	\psfrag{2}{\huge$2$}
	\psfrag{3}{\huge$3$}
	\psfrag{4}{\huge$4$}
	\psfrag{5}{\huge$\mathbf{5}$}
	\centerline{\scalebox{.5}{\includegraphics{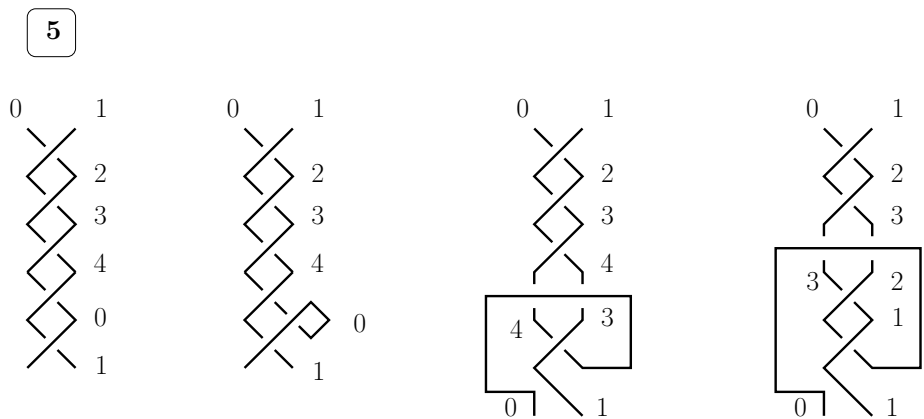}}}
	\caption{Decomposing Teneva's example.}\label{fig:decompteneva}
\end{figure}

Colorings of knots were introduced by Fox in \cite{CFox}. Minimum number of colors were introduced by Harary and the first author in \cite{Frank} where they conjectured that given a prime $p$ along with a reduced  alternating diagram of a knot with determinant $p$, then this diagram endowed with a non-trivial $p$-coloring would bear different colors on different arcs. This conjecture has now been proven true in \cite{msolis}. Still in \cite{Frank} it was remarked that an alternating knot of prime determinant $p$ could be non-trivially colored with less colors if one used diagrams other than alternating. This remark was made by Irina Teneva who provided the following example, see Figure \ref{fig:decompteneva}.
A few words about our figures seem to be in
order. Braid closure of each diagram in the figures is understood except where otherwise noted or clear from context. We do not depict these braid closures in order not to overburden the
figures. A circled number on the top left of a figure will stand for the modulus with respect to which the colorings are being considered.

Figure \ref{fig:decompteneva} shows a $T(2, 5)$ (whose determinant is $5$) endowed with a non-trivial $5$-coloring under four guises. The left-most diagram is a reduced alternating diagram thus complying with the Kauffman-Harary conjecture: distinct arcs bear distinct colors. But the rightmost diagram is a non-alternating diagram colored with less colors, $4$ to be precise. This is Teneva's example. The progression from the leftmost to the rightmost diagram is our decomposition of Teneva's example which allowed to generalize it to other $\sigma_1^n$'s, see \cite{kl}. We thus briefly remark that the so-called Teneva transformations are applicable to colored $\sigma_1^n$'s. Each Teneva transformation consists of a finite sequence of Reidemeister moves, each move accompanied by a corresponding local alteration in the colors in order to make the colorings before and after the Reidemeister move compliant with each other. A Teneva transformation begins then with a kink, it is followed by a first type III Reidemeister move which does not yet alter the number of colors, it is followed by a second type III Reidemeister move which removes one color, and this is repeated until the application of such moves introduces old colors instead of removing them. Figure \ref{fig:decompteneva} is the materialization of what we just described for $n=5$. For a more elaborate discussion see \cite{kl}.

\bigbreak

It should be clear that it is not possible to further apply Teneva transformations to the right-most diagram of Figure \ref{fig:decompteneva} in order to reduce the number of colors. On the other hand, for larger $n$'s, it is possible to effectively apply Teneva transformations to the $\sigma_1^m$'s left over from the application of the first set of Teneva transformations, or even from subsequent applications of Teneva transformations. We call this iterating the Teneva transformations. This is the Teneva Game.

Figure \ref{fig:t2-11} should help at this point. After a first application of the Teneva transformations, the diagram is split into two $\sigma_1^5$'s (Figure \ref{fig:t2-11}, second diagram from the left). Teneva transformations are then applied to these $\sigma_1^5$'s splitting them into $\sigma_1^2$'s (third diagram from the left).

\bigbreak

In the current article we iterate the Teneva transformations over torus knots of type $(2, p)$, for prime $p>7$ in order to obtain sharper upper bounds for the minimum number of colors of these $T(2, p)$ modulo $p$ (Theorem \ref{thm:thm}).  Furthermore, letting $mincol_n K$ stand for the minimum number of colors needed to realize a non-trivial $n$-coloring of knot/link $K$ modulo $n$, we show (Corollary \ref{cor:mincol11}) that
\[
mincol_{11}T(2, 11) = 5.
\]

\bigbreak

In \cite{lm} it was proven that given a non-split link $L$ admitting non-trivial $r$ colorings then $mincol_r L = mincol_p L$, for a certain prime $p | (r, \det L)$.
 For this reason the moduli with respect to which we will be working with here will be odd primes.

\bigbreak

This article is organized as follows. In Section \ref{sect:prelim} we introduce the objects which simplify the statement of our results, and we state these results in Theorem \ref{thm:thm}, and Corollaries \ref{cor:mincol11} and \ref{cor:mincol}. In Section \ref{sect:proofs} we prove the Theorem \ref{thm:thm} and Corollary \ref{cor:mincol11}. In Section \ref{sect:algo} we give an algorithm for a quick calculation of  the relevant parameters and prove Corollary \ref{cor:mincol}. In Section \ref{sect:estimates} we estimate the efficiency of our results for larger primes. In Section \ref{sect:preview} we offer the reader a preview of our next project.

\section{Terminology, notation, and main results.} \label{sect:prelim}

\noindent

In this Section we begin by introducing terminology and notation which will simplify the statement of the main results.

\begin{def.}[Lower Half]
Let $o$ be an odd integer. We call {\rm\bf lower half of $o$} the integer $k$ such that
\[
o = 2k+1.
\]
Our notation for this situation will be:
\[
lh (o) = k.
\]
For an even integer, the lower half of it coincides with the ordinary half.
\end{def.}

We remark that we could have also defined the {\rm\bf Upper Half} of an odd integer $o$ as the integer $k$ such that $o=2k-1$, together with the statement that for even integers it also coincides with the ordinary half, but we do not seem to need this object in the present article.

\begin{def.}[Sequence of Lower Halves]\label{def:slh}
Let $n$ be an integer greater than $4$. We define {\rm\bf Sequence of Lower Halves of $n$} as follows.

The first term of this sequence is:
\[
k_1 := lh(n).
\]

If $k_i \in \{ 2, 3, 4 \}$, then this is the last term of the sequence. Otherwise,
\[
k_{i+1} := lh(k_i).
\]

Roughly speaking, the {\rm Sequence of Lower Halves} is the sequence of iterates of the function $lh(\cdot )$ applied to $n$ and terminating at the first iterate in $\{ 2, 3, 4 \}$.

Our notation for the sequence of lower halves of $n$ is $LH(n)$, and $k_i$ for its terms. Given an integer $n$ along with its sequence of lower halves, $LH(n) = (k_1, k_2, \dots )$, we call {\rm\bf length of the sequence}, notation $l_n$, the number of terms in the sequence. We call {\rm\bf tail of the sequence}, notation $t_n$, the last term of the sequence, $k_{l_n}$.

For each prime $p>7$, we will use these quantities, $l_p$ and $t_p$, to form the upper bound in Theorem \ref{thm:thm}.
\end{def.}

Here are some examples:
\[
31 = 2\times 15 + 1 = 2\times [2\times 7 + 1] + 1 = 2\times [2\times (2\times 3 + 1) + 1] + 1
\]
\[
\text{LH $(31)$} = (15, 7, 3) \qquad \qquad \qquad \qquad l_{31} = 3  \qquad \qquad t_{31}= 3.
\]

\bigbreak

\[
37 = 2\times 18 + 1 = 2\times [2\times 9] + 1 = 2\times [2\times (2\times 4 + 1)] + 1
\]
\[
\text{LH $(37)$} = (18, 9, 4) \qquad \qquad l_{37} = 3  \qquad \qquad t_{37}= 4.
\]
\bigbreak

\[
41 = 2\times 20 + 1 = 2\times \{2\times 10 \} + 1  = 2\times \{2\times [2\times 5 ] \} + 1  = 2\times \{2\times [2\times (2\times 2 + 1) ] \} + 1
\]
\[
\text{LH $(41)$} = (20, 10, 5, 2) \qquad \qquad l_{41} = 4  \qquad \qquad t_{41}= 2.
\]

We remark that it would have been more natural to define the Sequences of Lower Halves terminating at $1$ but it does not seem to be helpful in our set up,
as will soon become clear.

\bigbreak

Here are the main results of this article. (In the sequel we retain the notation above concerning lower halves, sequences of lower halves, its terms, lengths of sequences and tails of sequences.)

\begin{thm}\label{thm:thm} Let $p>7$ be prime. Then,
\[
mincol_p T(2, p) \leq t_p + 2l_p -1,
\]
where $t_p$ and $l_p$ were introduced in Definition \ref{def:slh}.
\end{thm}

\bigbreak

The following result was announced in \cite{satoh}.

\begin{cor}\label{cor:mincol11}
\[
mincol_{11} T(2, 11) = 5.
\]
\end{cor}

\bigbreak

We remark that the following results are already known (see \cite{lm}, \cite{Oshiro}, and \cite{satoh}):

\[
mincol_3T(2, 3) = 3, \qquad \qquad \qquad \qquad mincol_5T(2, 5) = 4, \qquad \qquad \qquad \qquad mincol_7T(2, 7) = 4.
\]

\bigbreak

\begin{cor}\label{cor:mincol} For any prime $p > 7$,
\[
mincol_{p} T(2, p) \leq 2 \log_2 (p-1) - 1.
\]
\end{cor}

\bigbreak

\section{Proofs of Theorem \ref{thm:thm} and Corollary \ref{cor:mincol11}.} \label{sect:proofs}

\noindent

We prove Theorem \ref{thm:thm} in this Section. We start by considering some examples namely, $T(2, 11)$ and $T(2, 13)$. In the course of analyzing $T(2, 11)$ we prove Corollary \ref{cor:mincol11}.

\bigbreak

\begin{figure}[!ht]
	\psfrag{0}{\huge$0$}
	\psfrag{1}{\huge$1$}
	\psfrag{2}{\huge$2$}
	\psfrag{3}{\huge$3$}
	\psfrag{4}{\huge$4$}
	\psfrag{5}{\huge$5$}
	\psfrag{6}{\huge$6$}
	\psfrag{7}{\huge$7$}
	\psfrag{8}{\huge$8$}
	\psfrag{9}{\huge$9$}
	\psfrag{10}{\huge$10$}
	\psfrag{11}{\huge$\bf 11$}
	\psfrag{7-10 removed}{\huge$7 \text{ through } 10 \text{ removed}$}
	\psfrag{k1-1 colors}{\huge$\text{i.e., } k_1-1 \text{ colors removed}$}
	\psfrag{4 removed}{\huge$4 \text{ removed}$}
	\psfrag{k2-1 colors}{\huge$\text{i.e., } k_2-1 \text{ colors removed}$}
	\psfrag{5 removed}{\huge$5 \text{ removed}$}
	\psfrag{1 colors}{\huge$\text{i.e., } 1 \text{ color removed}$}
	\psfrag{k1=5}{\huge$k_1=5$}
	\psfrag{k2=2}{\huge$k_2=2$}
	\psfrag{b2}{\huge$b_2$}
	\psfrag{b3}{\huge$b_3$}
	\psfrag{b4}{\huge$b_4$}
	\psfrag{b5}{\huge$b_5$}
	\psfrag{b6}{\huge$b_6$}
	\centerline{\scalebox{.45}{\includegraphics{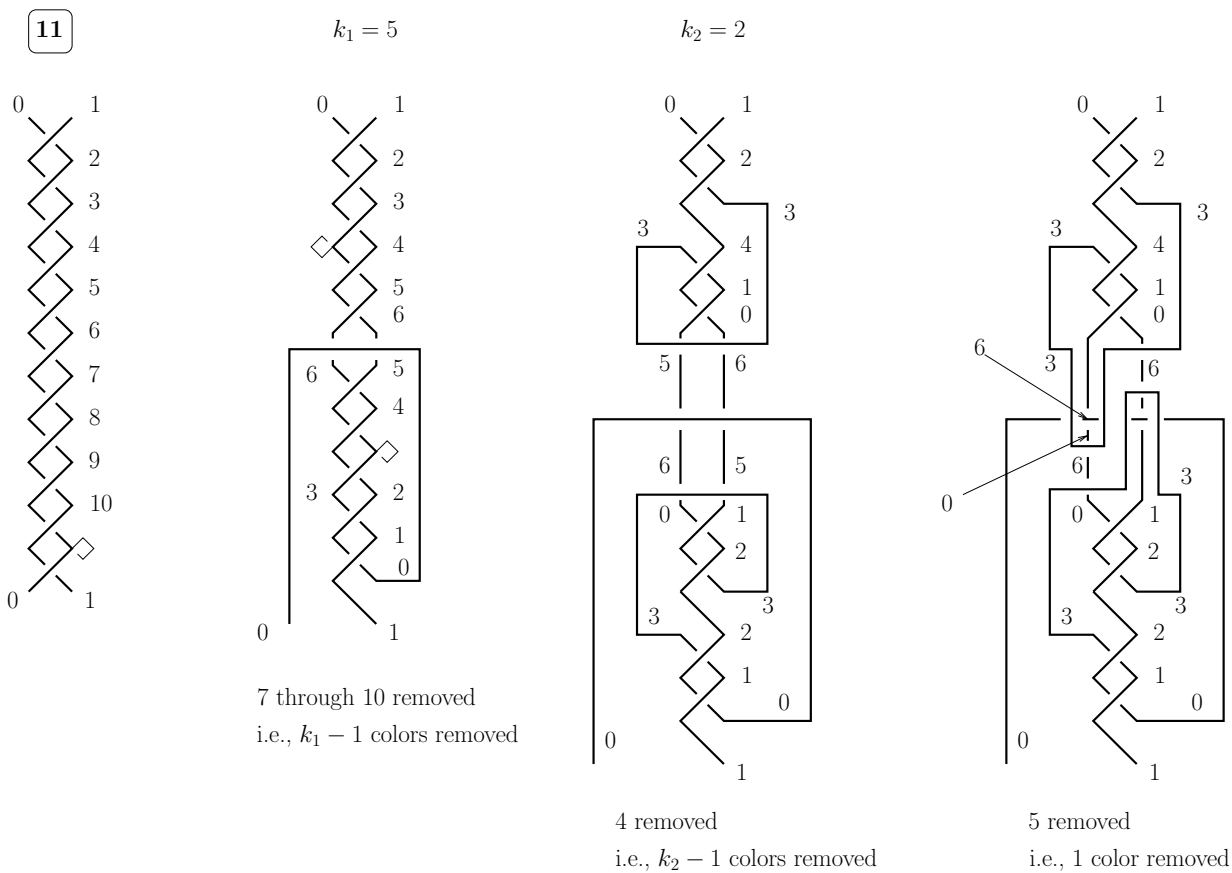}}}
	\caption{$T(2, 11)$ admits a non-trivial $11$-coloring with $5$ colors:  $0, 1, 2, 3, 6$ (rightmost diagram). This is the minimum number of colors
for this knot.}\label{fig:t2-11}
\end{figure}

We remark that
\[
LH(11) = (5, 2)  \qquad \qquad \qquad \qquad k_1 = 5 \qquad k_2 = 2 \qquad \qquad \qquad \qquad l_{11}=2 \qquad t_{11}=2.
\]
Consider Figure \ref{fig:t2-11}. (Here and in the sequel, kinks drawn in a thinner line indicate the beginning of the Teneva transformation depicted in the diagram to the right of the one where the kink is found.) The first set of Teneva transformations performed on $T(2, 11)$ takes an arc over the diagram, splitting it into two $\sigma_1^5$'s. This is depicted in the second diagram from the left in Figure \ref{fig:t2-11}. This transformation removes colors $7$ through $10$; this is the removal of $k_1 - 1$ colors for $k_1 = 5$. Why did this first set of Teneva transformations take this arc to split this diagram like it did? Well, had this arc gone over more crossings or less crossings than it did then it would have eliminated less colors. In the general case of a $\sigma_1^{2k+1}$ endowed with a non-trivial $(2k+1)$-coloring, the maximum reduction of colors by one set of Teneva transformations is accomplished when the diagram is split into two $\sigma_1^k$'s, yielding $k+2$ distinct colors (here k=5). This is carefully discussed in \cite{kl}.

Let us resume the analysis of Figure \ref{fig:t2-11}, on the second diagram from the left. The diagram is split into two $\sigma_1^5$'s. Teneva transformations can be applied effectively to each of these $\sigma_1^5$'s, reducing their common set of colors. The starting points for the new Teneva transformations are indicated by the kinks and the result is depicted on the third diagram from the left. Color $4$ has been removed from this diagram; this is the removal of $k_2-1$ colors for $k_2=2$. Finally, the rightmost diagram depicts a transformation which removes color $5$ leaving the diagram with $5$ colors: $0, 1, 2, 3, 6$. According to \cite{lm}, Theorem $1.4$, $(4)$ (see also \cite{Saito}), this is the least number of colors such a knot can exhibit in a non-trivial $11$-coloring. Therefore
\[
mincol_{11} T(2, 11) = 5,
\]
which concludes the proof of Corollary \ref{cor:mincol11}.

\bigbreak

\begin{figure}[!ht]
	\psfrag{0}{\huge$0$}
	\psfrag{1}{\huge$1$}
	\psfrag{2}{\huge$2$}
	\psfrag{3}{\huge$3$}
	\psfrag{4}{\huge$4$}
	\psfrag{5}{\huge$5$}
	\psfrag{6}{\huge$6$}
	\psfrag{7}{\huge$7$}
	\psfrag{8}{\huge$8$}
	\psfrag{9}{\huge$9$}
	\psfrag{10}{\huge$10$}
	\psfrag{11}{\huge$11$}
	\psfrag{12}{\huge$12$}
	\psfrag{13}{\huge$\bf 13$}
	\psfrag{8-12 removed}{\huge$8 \text{ through } 12 \text{ removed}$}
	\psfrag{k1-1 colors}{\huge$\text{i.e., } k_1-1 \text{ colors removed}$}
	\psfrag{5 removed}{\huge$5 \text{ removed}$}
	\psfrag{k2-2 colors}{\huge$\text{i.e., } k_2-2 \text{ colors removed}$}
	\psfrag{5 removed}{\huge$5 \text{ removed}$}
	\psfrag{1 colors}{\huge$\text{i.e., } 1 \text{ color removed}$}
	\psfrag{k1=6}{\huge$k_1=6$}
	\psfrag{k2=3}{\huge$k_2=3$}
	\psfrag{b2}{\huge$b_2$}
	\psfrag{b3}{\huge$b_3$}
	\psfrag{b4}{\huge$b_4$}
	\psfrag{b5}{\huge$b_5$}
	\psfrag{b6}{\huge$b_6$}
	\centerline{\scalebox{.5}{\includegraphics{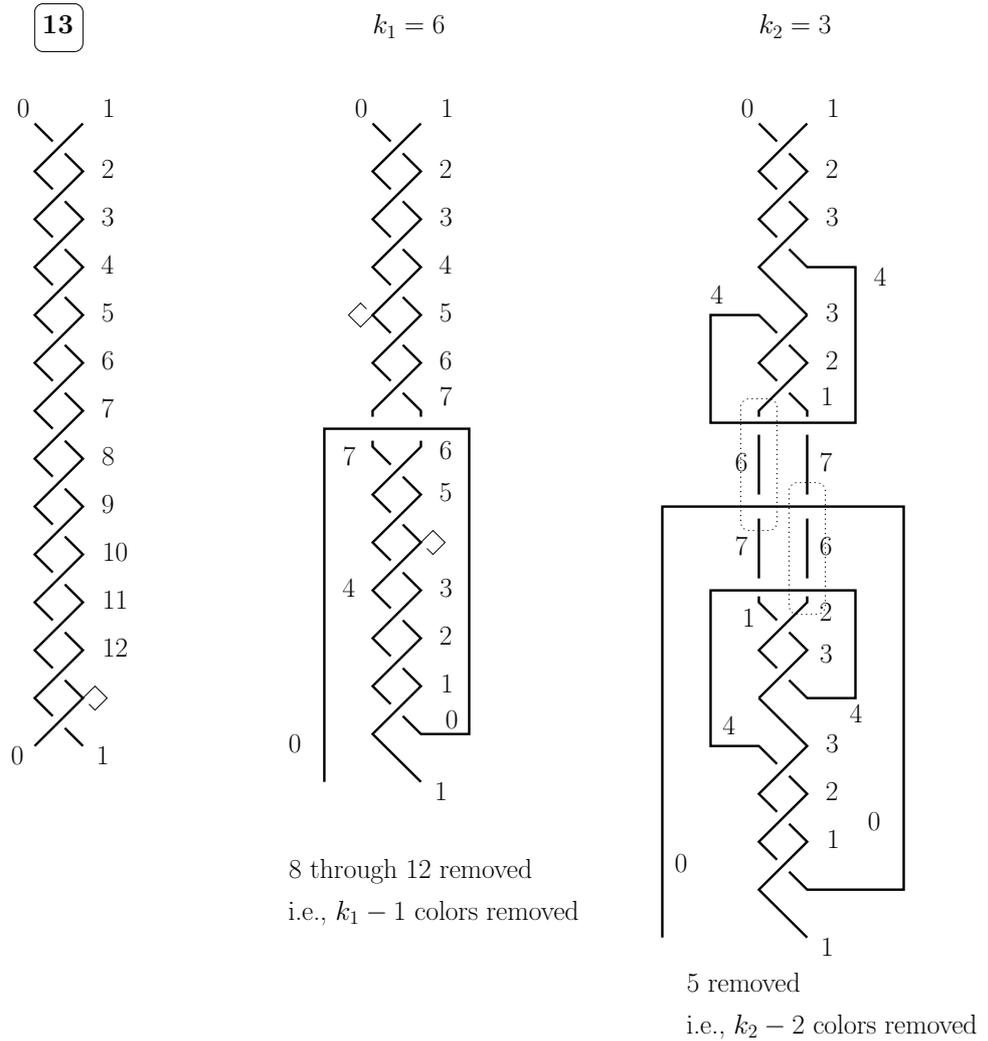}}}
	\caption{Towards realizing a non-trivial $13$-coloring of $T(2, 13)$ with $6$ colors (one color short, here).}\label{fig:t2-13}
\end{figure}

The following features of this example should be retained. With the indicated values of $k_1$ and $k_2$, the first set of Teneva transformations (the passage from the leftmost diagram to the next diagram in Figure \ref{fig:t2-11}) removes $k_1-1$ colors namely, colors $k_{1}+2$ through $2k_1$. The next set of Teneva transformations (the passage from the second diagram from the left to the third) removes $k_2 -1$ colors, namely colors $k_2 +2$ through $2k_2$ (which in this case is just $1$ color). It should also be noted that $k_1$ is odd; were it even then this second set of Teneva transformations would have removed only $k_2 -2$ as we will note in the next example, see Figure \ref{fig:t2-13}. Finally, the passage from the third from the left diagram to the right-most diagram removes color $k_1$. The total number of colors in the end is $t_{11} + 2l_{11}-1=5$

\bigbreak

Let us now analyze Figure \ref{fig:t2-13}. Here we have
\[
LH(13) = (6, 3) \qquad \qquad \qquad \qquad k_1=6 \quad k_2 = 3 \qquad \qquad \qquad \qquad l_{13}=2 \qquad t_{13}=3.
\]

\begin{figure}[!ht]
	\psfrag{0}{\huge$0$}
	\psfrag{1}{\huge$1$}
	\psfrag{2}{\huge$2$}
	\psfrag{4}{\huge$4$}
	\psfrag{6}{\huge$6$}
	\psfrag{7}{\huge$7$}
	\psfrag{8}{\huge$8$}
	\psfrag{13}{\huge$\bf 13$}
	\centerline{\scalebox{.5}{\includegraphics{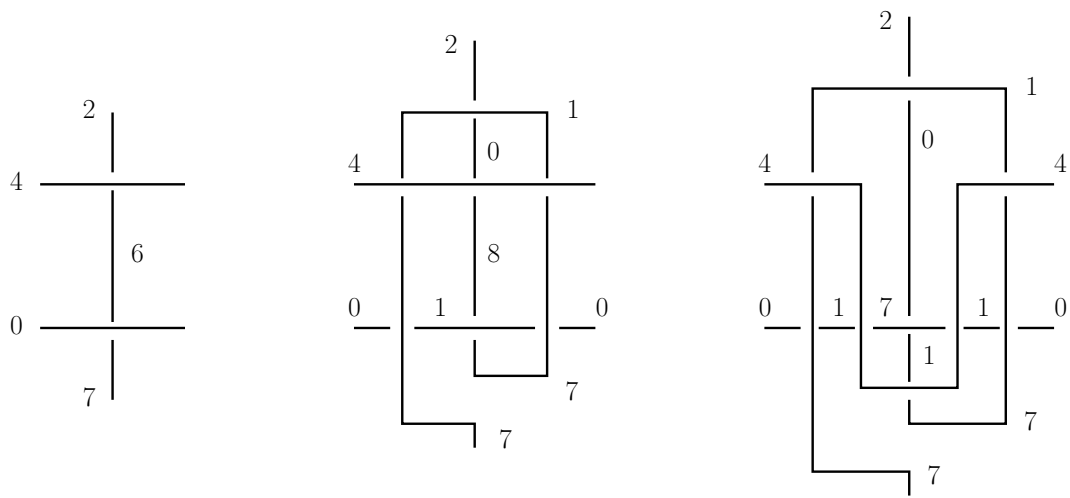}}}
	\caption{Boxed regions in Figure \ref{fig:t2-13}: manoeuvering to eliminate $6$.}\label{fig:t2-13localtransf}
\end{figure}

\begin{figure}[!ht]
	\psfrag{0}{\huge$0$}
	\psfrag{1}{\huge$1$}
	\psfrag{2}{\huge$2$}
	\psfrag{3}{\huge$3$}
	\psfrag{4}{\huge$4$}
	\psfrag{5}{\huge$5$}
	\psfrag{6}{\huge$6$}
	\psfrag{7}{\huge$7$}
	\psfrag{8}{\huge$8$}
	\psfrag{9}{\huge$9$}
	\psfrag{10}{\huge$10$}
	\psfrag{11}{\huge$11$}
	\psfrag{12}{\huge$12$}
	\psfrag{13}{\huge$\bf 13$}
	\centerline{\scalebox{.5}{\includegraphics{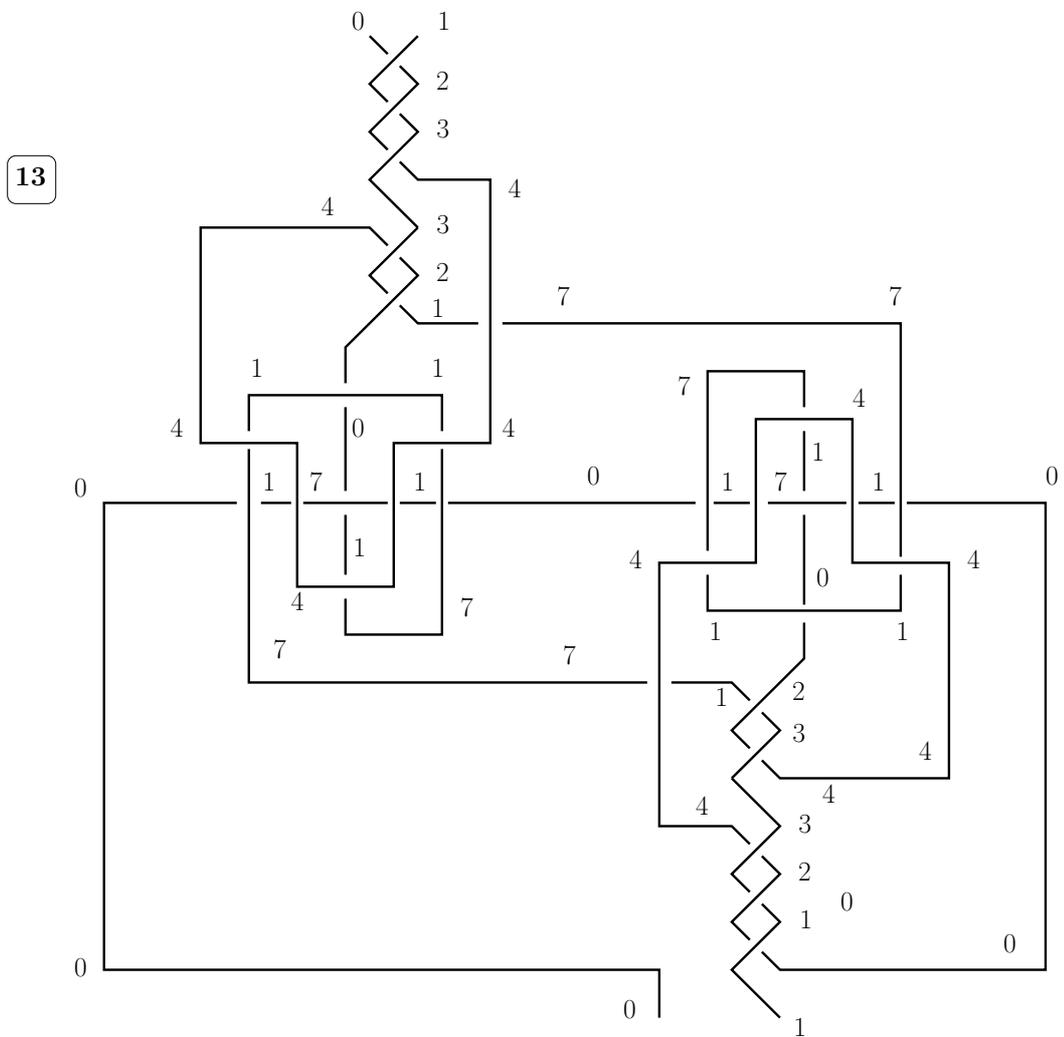}}}
	\caption{Introducing the treated boxed regions into the left-most diagram of Figure \ref{fig:t2-13}: $T(2, 13)$ admits a non-trivial $13$-coloring with $6$ colors.}\label{fig:t2-13final}
\end{figure}

The relevant difference from the preceding case is that here $k_1$ is even. Although the first set of Teneva transformations (the passage from the left-most to the center diagram) involve again the removal of $k_1-1$ colors namely, colors $k_1+2$ through $2k_1$, the second set of Teneva transformations involve the removal of $k_2-2$ colors namely, colors $k_2+2$ through $2k_2-1$ (which in this case is only $1$ color). In general, the following occurs. Suppose $k_i$ is not the tail of the sequence of lower halves of the modulus at issue. If $k_i$ is odd then the $(i+1)$-th set of Teneva transformations will remove $k_{i+1}-1$ colors; otherwise it will remove $k_{i+1}-2$ colors. Back to Figure \ref{fig:t2-13}, the boxed regions in the right-most diagram indicate where to operate in order to remove color $k_1$ which in this case is $6$. Since this is graphically more involved then the removal of the odd $k_1$ (see Figure \ref{fig:t2-11}), it is displayed in Figure \ref{fig:t2-13localtransf}. Inspiration for these moves comes from $\cite{Oshiro}$. Figure \ref{fig:t2-13final} incorporates this removal of color $6$, broken down in Figure \ref{fig:t2-13localtransf}, thus displaying the final result. The total number of colors in the end is $t_{13}+2l_{13}-1 = 6$.

\bigbreak

We are now ready for the general case. We remark that removal of color $k_1$ in the general case is obtained in a straightforward way from the two preceding examples. For instance, for even $k_1$ with Figure \ref{fig:t2-13localtransf} in mind, $6\rightarrow k_1$, $7\rightarrow k_1+1$, and $4\rightarrow \frac{k_1}{2}+1$. In order to
prove Theorem \ref{thm:thm} we use the following Lemma.

Only Teneva transformations are considered in Lemma \ref{lem:lem}. In particular, color $k_1$ is accounted for. Color $k_1$ can be removed by the procedures described above after all the sets of Teneva transformations have been applied.

\bigbreak

\begin{lem}\label{lem:lem}
Let $p$ be a prime greater than $7$ with
\[
LH(p) = (k_1, k_2, \dots , k_{i-1}, k_i, \dots , t_p) .
\]
For $i= 1, 2, \dots , l_p$, after applying the $i$-th set of Teneva Transformations as described above, we are left with $k_i+2i$ colors namely,
\[
0, 1, 2, \dots , k_{i}-1, \qquad k_i, k_{i}+1, \qquad k_{i-1}, k_{i-1}+1, \qquad \dots \quad , \quad k_{2}, k_{2}+1, \qquad  k_{1}, k_{1}+1
\]
\end{lem}
Proof. By induction on $i$. The proof for $i=1$ is essentially contained in the discussion of Figure \ref{fig:t2-11} and is carefully discussed in \cite{kl}.

Consider $1\leq i < l_p$ such that the colors left after the application of the $i$-th set of Teneva transformations are
\[
0, 1, 2, \dots , k_{i}-1, \qquad k_i, k_{i}+1, \qquad k_{i-1}, k_{i-1}+1, \qquad \dots \quad , \quad k_{2}, k_{2}+1,  \qquad k_{1}, k_{1}+1
\]
i.e., $k_i+2i$ colors.
Colors $1, 2, \dots , k_{i}-1, k_i, k_{i}+1$  show up in $\sigma_1^{k_i-1}$'s as shown in Figures \ref{fig:sigmakiodd} and \ref{fig:sigmakieven} (their left-most diagrams only, for now) and nowhere else. 
\begin{enumerate}
\item Assume $k_i$ is odd. The application of the $(i+1)$-th set of Teneva transformations to these $\sigma_1^{k_i-1}$'s (see Figure \ref{fig:sigmakiodd}) removes colors $k_{i+1}+2$ through $2k_{i+1}(=k_i-1)$ leaving behind colors
\[
0, 1, \dots , k_{i+1}-1, \qquad k_{i+1}, k_{i+1}+1, \qquad k_i, k_i+1,  \quad \dots \quad , \quad  k_{2}, k_{2}+1, \quad k_{1}, k_{1}+1
\]
i.e., $k_{i+1}+2(i+1)$ colors.
\begin{figure}[!ht]
	\psfrag{0}{\huge$0$}
	\psfrag{1}{\huge$1$}
	\psfrag{2}{\huge$2$}
	\psfrag{3}{\huge$3$}
	\psfrag{4}{\huge$4$}
	\psfrag{5}{\huge$5$}
	\psfrag{6}{\huge$6$}
	\psfrag{7}{\huge$7$}
	\psfrag{8}{\huge$8$}
	\psfrag{9}{\huge$9$}
	\psfrag{10}{\huge$10$}
	\psfrag{11}{\huge$11$}
	\psfrag{12}{\huge$12$}
	\psfrag{13}{\huge$\bf 13$}
	\psfrag{8-12 removed}{\huge$8 \text{ through } 12 \text{ removed}$}
	\psfrag{k1-1 colors}{\huge$\text{i.e., } k_1-1 \text{ colors removed}$}
	\psfrag{5 removed}{\huge$5 \text{ removed}$}
	\psfrag{k2-2 colors}{\huge$\text{i.e., } k_2-2 \text{ colors removed}$}
	\psfrag{5 removed}{\huge$5 \text{ removed}$}
	\psfrag{1 colors}{\huge$\text{i.e., } 1 \text{ color removed}$}
	\psfrag{ki=2ki+1+1}{\huge$\mathbf{k_i=2k_{i+1}+1}$}
	\psfrag{ki +1}{\huge$k_i+1$}
	\psfrag{ki+1}{\huge$k_{i+1}$}
	\psfrag{ki+1+1}{\huge$k_{i+1}+1$}
	\psfrag{ki+1+2}{\huge$k_{i+1}+2$}
	\psfrag{ki+1-1}{\huge$k_{i+1}-1$}
	\psfrag{ki+1-2}{\huge$k_{i+1}-2$}
	\psfrag{...}{\huge$\cdots$}
	\psfrag{ki}{\huge$k_i$}
	\psfrag{ki-1}{\huge$k_{i}-1$}
	\psfrag{ki-2}{\huge$k_{i}-2$}
	\psfrag{k2=3}{\huge$k_2=3$}
	\psfrag{ki+1+2-ki-1 removed}{\huge$k_{i+1}+2 \text{ through } k_{i}-1 \text{ removed}$}
	\psfrag{b2}{\huge$b_2$}
	\psfrag{b3}{\huge$b_3$}
	\psfrag{b4}{\huge$b_4$}
	\psfrag{b5}{\huge$b_5$}
	\psfrag{b6}{\huge$b_6$}
	\centerline{\scalebox{.45}{\includegraphics{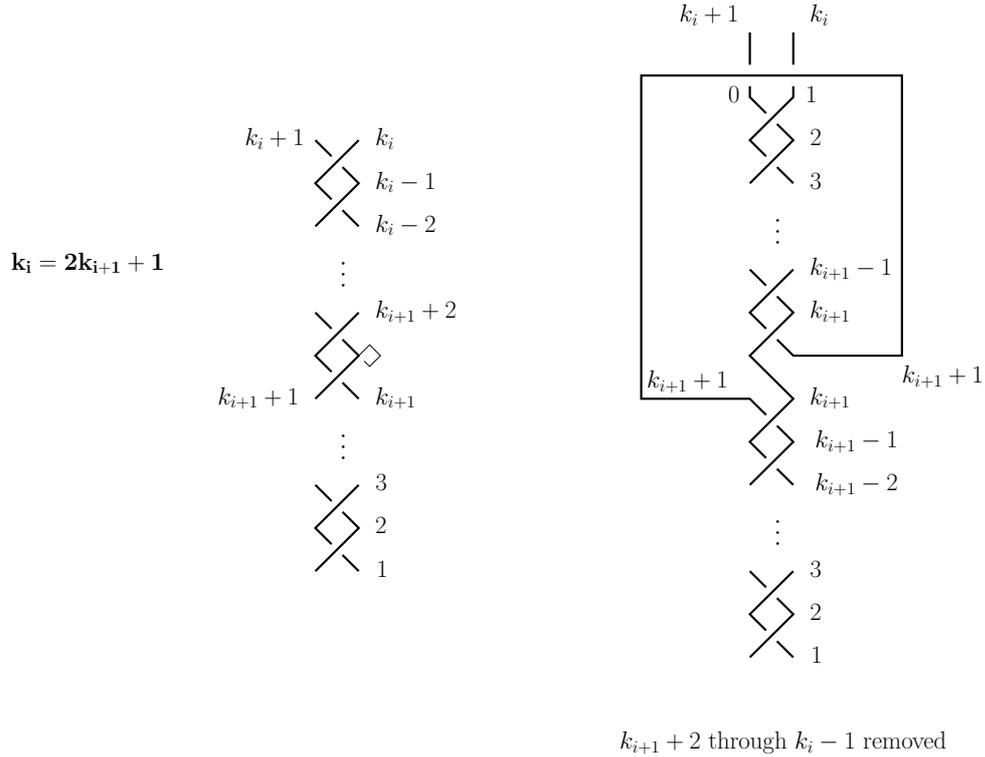}}}
	\caption{The $(i+1)$-th Teneva transformation on a colored $\sigma_1^{k_i-1}$ for odd $k_i$ removes $k_{i+1}-1$ colors.}\label{fig:sigmakiodd}
\end{figure}
\item Assume $k_i$ is even. The application of the $(i+1)$-th set of Teneva transformations to these $\sigma_1^{k_i-1}$'s (see Figure \ref{fig:sigmakieven}) removes colors $k_{i+1}+2$ through $2k_{i+1}-1(=k_i-1)$ leaving behind colors
\[
0, 1, \dots , k_{i+1}-1, \qquad k_{i+1}, k_{i+1}+1, \qquad k_i, k_i+1,  \quad \dots \quad , \quad  k_{2}, k_{2}+1, \quad k_{1}, k_{1}+1
\]
i.e., $k_{i+1}+2(i+1)$ colors.
\begin{figure}[!ht]
	\psfrag{0}{\huge$0$}
	\psfrag{1}{\huge$1$}
	\psfrag{2}{\huge$2$}
	\psfrag{3}{\huge$3$}
	\psfrag{4}{\huge$4$}
	\psfrag{5}{\huge$5$}
	\psfrag{6}{\huge$6$}
	\psfrag{7}{\huge$7$}
	\psfrag{8}{\huge$8$}
	\psfrag{9}{\huge$9$}
	\psfrag{10}{\huge$10$}
	\psfrag{11}{\huge$11$}
	\psfrag{12}{\huge$12$}
	\psfrag{13}{\huge$\bf 13$}
	\psfrag{8-12 removed}{\huge$8-12 \text{ removed}$}
	\psfrag{k1-1 colors}{\huge$\text{i.e., } k_1-1 \text{ colors removed}$}
	\psfrag{5 removed}{\huge$5 \text{ removed}$}
	\psfrag{k2-2 colors}{\huge$\text{i.e., } k_2-2 \text{ colors removed}$}
	\psfrag{5 removed}{\huge$5 \text{ removed}$}
	\psfrag{1 colors}{\huge$\text{i.e., } 1 \text{ color removed}$}
	\psfrag{ki=2ki+1}{\huge$\mathbf{k_i=2k_{i+1}}$}
	\psfrag{ki +1}{\huge$k_i+1$}
	\psfrag{ki+1}{\huge$k_{i+1}$}
	\psfrag{ki+1+1}{\huge$k_{i+1}+1$}
	\psfrag{ki+1+2}{\huge$k_{i+1}+2$}
	\psfrag{ki+1-1}{\huge$k_{i+1}-1$}
	\psfrag{ki+1-2}{\huge$k_{i+1}-2$}
	\psfrag{...}{\huge$\cdots$}
	\psfrag{ki}{\huge$k_i$}
	\psfrag{ki-1}{\huge$k_{i}-1$}
	\psfrag{ki-2}{\huge$k_{i}-2$}
	\psfrag{k2=3}{\huge$k_2=3$}
	\psfrag{ki+1+2-ki-1 removed}{\huge$k_{i+1}+2 \text{ through } k_{i}-1 \text{ removed}$}
	\psfrag{b2}{\huge$b_2$}
	\psfrag{b3}{\huge$b_3$}
	\psfrag{b4}{\huge$b_4$}
	\psfrag{b5}{\huge$b_5$}
	\psfrag{b6}{\huge$b_6$}
	\centerline{\scalebox{.45}{\includegraphics{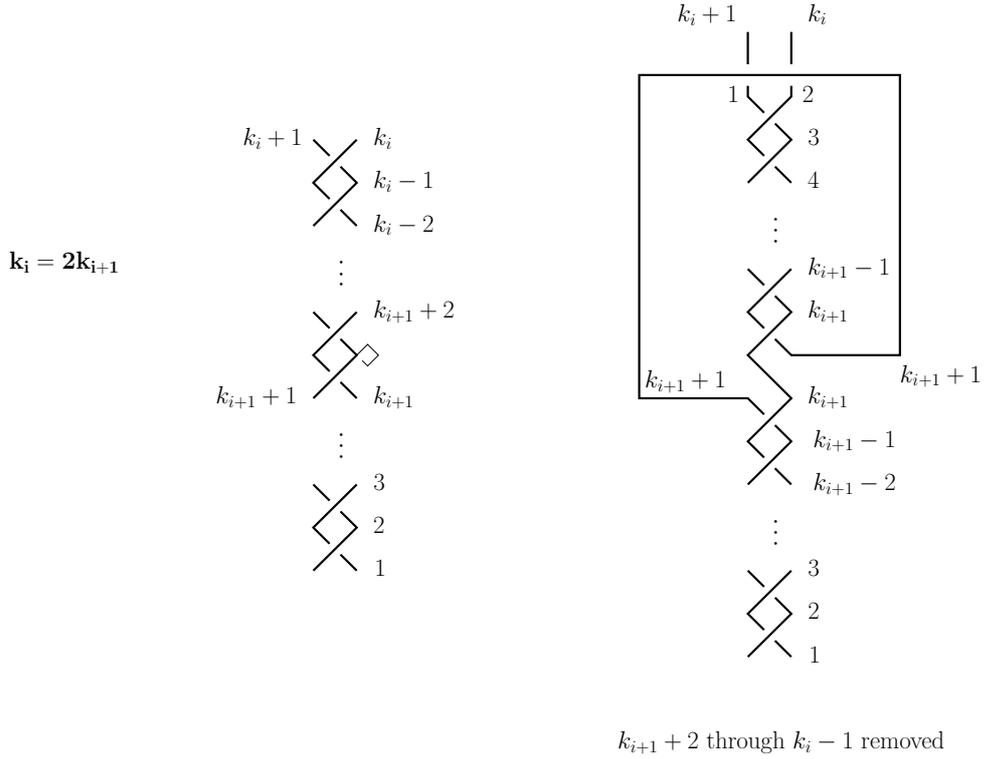}}}
	\caption{The $(i+1)$-th Teneva transformation on a colored $\sigma_1^{k_i-1}$ for even $k_i$ removes $k_{i+1}-2$ colors.}\label{fig:sigmakieven}
\end{figure}
\end{enumerate}

This concludes the proof of  Lemma \ref{lem:lem}.

$\hfill \blacksquare$

\bigbreak

Proof of Theorem \ref{thm:thm}.

Leaning on Lemma \ref{lem:lem}, given a prime $p>7$, after applying $l_p$ sets of Teneva transformations to $T(2, p)$ as described above, we are left with $t_p+2l_p$ colors. Finally we remove color $k_{1}$ by the procedures described before. We are then left with $t_p+2l_p-1$ colors. This concludes the proof of Theorem \ref{thm:thm}.

$\hfill \blacksquare$

\bigbreak

\section{An algorithm for calculating the length and the tail of sequences of lower halves.} \label{sect:algo}

\noindent

\begin{prop}\label{prop:algo}

Consider the positive odd integer $n=2k+1$ along with its $2$-adic expansion:
\[
n= 2^{e_1} + 2^{e_2} + \cdots + 2^{e_{N_n}} + 1 \qquad \qquad \text{ with } e_1 > e_2 > \cdots > e_{N_n} > 0
\]

Then
\begin{enumerate}
\item If $e_1 - e_2 = 1$ then
\[
l_n = e_1-1 \qquad \qquad \text{ and } \qquad \qquad t_n=3
\]

\item If $e_1 - e_2 = 2$ then
\[
l_n = e_1-1 \qquad \qquad  \text{ and }  \qquad \qquad  t_n=2
\]

\item If $e_1 - e_2 > 2$ or $N_n=1$, then
\[
l_n = e_1-2 \qquad \qquad  \text{ and }  \qquad \qquad t_n=4
\]
\end{enumerate}
\end{prop}Proof: We keep the notation of the statement of the Proposition and rewrite the $2$-diadic expansion of $n$ in the following way.
\[
n =2^{s_{1}+s_{2}+ \cdots + s_{N_n-1}+s_{N_n}}+ 2^{s_{1}+ \cdots + s_{N_n-1}+s_{N_n}} + \cdots  + 2^{s_{N_n-1}+s_{N_n}} + 2^{s_{N_n}} + 1
\]
where
\[
s_{N_n}:=e_{N_n} \qquad s_{N_n-1}:=e_{N_n-1}-e_{N_n} \qquad \cdots \qquad s_2:=e_2-e_3 \qquad s_1:=e_1-e_2
\]
\begin{enumerate}
\item If $e_1 - e_2 = 1$ then
\begin{align*}
n &= 2^{1+e_2} +  2^{e_2} +  2^{e_3} + \cdots +  2^{e_{N_n}} + 1 = \\
&= 2^{1+s_2+s_3+\cdots + s_{N_n}} + 2^{s_2+s_3+\cdots + s_{N_n}} + 2^{s_3+\cdots + s_{N_n}} + \cdots + 2^{s_{N_n}} + 1 =  \\
&= 2^{1+s_2+s_3+\cdots + s_{N_n}} + 2^{s_2+s_3+\cdots + s_{N_n}} + 2^{s_3+\cdots + s_{N_n}} + \cdots + 2^{s_{N_n}} + 1 = \\
&= 2^{s_{N_n}}\big( 2^{1+s_2+s_3+\cdots + s_{{N_n}-1}} + 2^{s_2+s_3+\cdots + s_{{N_n}-1}} + 2^{s_3+\cdots + s_{{N_n}-1}} + \cdots +  1\big) + 1 =\\
&= 2^{s_{N_n}}\bigg( 2^{s_{{N_n}-1}}\big(2^{1+s_2+s_3+\cdots + s_{{N_n}-2}} + 2^{s_2+s_3+\cdots + s_{{N_n}-2}} + 2^{s_3+\cdots + s_{{N_n}-2}} + \cdots + 1\big)  + 1\bigg) + 1 = \\
&= \cdots = 2^{s_{N_n}}\Bigg( 2^{s_{{N_n}-1}}\Big( \cdots 2^{s_2}\big( 2^1 + 1\big) \cdots  + 1\Big)  + 1\Bigg) + 1 =
\end{align*}
So
\[
l_n = s_{N_n} + s_{{N_n}-1} + \cdots + s_2 = e_2 = e_1 -1 \qquad \qquad \qquad \qquad t_n = 3
\]

\item If $e_1 - e_2 = 2$ then
\begin{align*}
n &= 2^{2+e_2} +  2^{e_2} +  2^{e_3} + \cdots +  2^{e_{N_n}} + 1 = \\
&= 2^{2+s_2+s_3+\cdots + s_{N_n}} + 2^{s_2+s_3+\cdots + s_{N_n}} + 2^{s_3+\cdots + s_{N_n}} + \cdots + 2^{s_{N_n}} + 1 =  \\
&= \cdots = 2^{s_{N_n}}\Bigg( 2^{s_{{N_n}-1}}\Big( \cdots 2^{s_2}\big( 2\times 2 + 1\big) \cdots  + 1\Big)  + 1\Bigg) + 1 =
\end{align*}
So
\[
l_n = s_{N_n} + s_{{N_n}-1} + \cdots + s_2 = e_2 = e_1 -1 \qquad \qquad \qquad \qquad t_n = 2
\]

\item We leave the $N_n=1$ instance for the reader. If $e_1 - e_2 > 2$ then
\begin{align*}
n &= 2^{e_1} +  2^{e_2} +  2^{e_3} + \cdots +  2^{e_{N_n}} + 1 = \\
&= 2^{s_1+s_2+s_3+\cdots + s_{N_n}} + 2^{s_2+s_3+\cdots + s_{N_n}} + 2^{s_3+\cdots + s_{N_n}} + \cdots + 2^{s_{N_n}} + 1 =  \\
&= \cdots = 2^{s_{N_n}}\Bigg( 2^{s_{{N_n}-1}}\Big( \cdots 2^{s_2}\big( 2\times 2^{s_1-1} + 1\big) \cdots  + 1\Big)  + 1\Bigg) + 1 =
\end{align*}
So
\[
l_n = s_{N_n} + s_{{N_n}-1} + \cdots + s_2 + s_{1}-2= e_1 -2 \qquad \qquad \qquad \qquad t_n = 4
\]
\end{enumerate}
$\hfill \blacksquare$

With Proposition \ref{prop:algo} we give the following algorithm for retrieving $t_p$ and $l_p$, given a prime $p>7$. Once the two largest exponents in the $2$-diadic expansion of $p$ are obtained, the absolute value of their difference along with the largest exponent, say $e_1$, tell us right away the parameters we need. Specifically, if the indicated difference is $1$, then the tail is $3$ and the length is $e_1-1$; if the difference is $2$ then the tail is $2$ and the length is $e_1-1$; finally if the difference is greater than $2$, then the tail is $4$ and the length is $e_1-2$.

\bigbreak

\noindent
{\bf Proof of Corollary \ref{cor:mincol}:}
\bigbreak
Except for the replacement of odd integer $n$ for prime $p$ (greater than $7$), we keep the notation of Proposition \ref{prop:algo}.

We will prove that
\[
2^{\frac{t_p+2l_p}{2}} \leq p-1
\]
Together with Theorem \ref{thm:thm} this will prove the statement of Corollary \ref{cor:mincol}.

We will now consider each of the three instances in the statement of Proposition \ref{prop:algo}.

\begin{enumerate}
\item If $e_1 - e_2 = 1$ then
\[
l_p = e_1-1 \qquad \qquad \text{ and } \qquad \qquad t_p=3
\]

which implies
\[
2^{\frac{t_p+2l_p}{2}} = 2^{\frac{3+2(e_1-1)}{2}}   = 2^{\frac{1}{2} + e_1} = \sqrt{2}\times 2^{e_1} \leq \bigg(1+\frac{1}{2}\bigg)\times 2^{e_1} = 2^{e_1} + 2^{e_1-1} \leq p-1
\]

\item If $e_1 - e_2 = 2$ then
\[
l_p = e_1-1 \qquad \qquad  \text{ and }  \qquad \qquad  t_p=2
\]

which implies
\[
2^{\frac{t_p+2l_p}{2}} = 2^{\frac{2+2(e_1-1)}{2}} =  2^{e_1} \leq p-1
\]

\item If $e_1 - e_2 > 2$ or $N_p=1$, then
\[
l_p = e_1-2 \qquad \qquad  \text{ and }  \qquad \qquad t_p=4
\]

which implies
\[
2^{\frac{t_p+2l_p}{2}} = 2^{\frac{4+2(e_1-2)}{2}} =  2^{e_1} \leq p-1
\]
\end{enumerate}
This concludes the proof of Corollary \ref{cor:mincol}.
$\hfill \blacksquare$
\section{Other estimates} \label{sect:estimates}

\noindent

In this Section we estimate the ratio between the number of colors obtained from the application of our procedure and the total number of colors available, showing that this ratio tends to decrease with increasing number of colors available.

We fix a modulus, an odd prime $p=2k+1$. We let $e_1$ and $e_2$ be integers such that
\[
2^{e_1} < p < 2^{e_1+1} \qquad \qquad \qquad \qquad 2^{e_2} < p - 2^{e_1} < 2^{e_2+1}
\]

The reduced number of colors is $t_p+2l_p - 1$, as proved in Theorem \ref{thm:thm}.

\begin{itemize}
\item If $e_1-e_2= 1$ then $l_p=e_1-1$ and $t_p=3$. Then, $t_p+2l_p-1 = 2+2l_p$ and

\[
\frac{1}{2}\frac{1+l_p}{2^{l_p}} = \frac{2+2l_p}{2^{2+l_p}} <  \frac{t_p+2l_p-1}{p} = \frac{2+2l_p}{p} < \frac{2+2l_p}{2^{1+l_p}}= \frac{1+l_p}{2^{l_p}}
\]

\item If $e_1-e_2= 2$, then $l_p=e_1-1$ and $t_p=2$. Then, $t_p+2l_p-1 = 1+2l_p$ and

\[
\frac{1}{2}\frac{1/2+l_p}{2^{l_p}} = \frac{1+2l_p}{2^{2+l_p}} <  \frac{t_p+2l_p-1}{p} = \frac{1+2l_p}{p} < \frac{1+2l_p}{2^{1+l_p}}= \frac{1/2+l_p}{2^{l_p}} < \frac{1+l_p}{2^{l_p}}
\]

\item If $e_1-e_2 > 2$ or $N_n=1$, then $l_p=e_1-2$ and $t_p=4$. Then, $t_p+2l_p-1 = 3+2l_p$ and

\[
\frac{1}{2}\frac{3/2+l_p}{2^{1+l_p}} = \frac{3+2l_p}{2^{3+l_p}} <  \frac{t_p+2l_p-1}{p} = \frac{3+2l_p}{p} < \frac{3+2l_p}{2^{2+l_p}}= \frac{1}{2} \frac{3/2+l_p}{2^{l_p}} < \frac{1}{2} \frac{2+2l_p}{2^{l_p}}  = \frac{1+l_p}{2^{l_p}}
\]

\end{itemize}

Consider then the  function
\[
f(x) = \frac{1+x}{2^x}
\]
for $x>1$.

We have
\[
f'(x) =\frac{2^x - (1+x)2^x\ln 2}{(2^x)^2} = \frac{1-(1+x)\ln 2}{2^x} < 0 , \quad \text{ for $x>1$}
\]

Hence $f(l_p)$ decreases with $l_p$ and
\[
\frac{t_p + 2l_p - 1}{p} < f(l_p).
\]

\bigbreak

Here is some experimenting with actual figures.
\[
p=11 \qquad \qquad \qquad LH(11)=(5, 2) \qquad \qquad \qquad \frac{2+2\times 2 -1}{11} < 0.455
\]

\[
p=13 \qquad \qquad \qquad LH(13)=(6, 3) \qquad \qquad \qquad \frac{3+2\times 2 -1}{13} < 0.462
\]

\[
p=331 \qquad \quad   LH(331) = (165, 82, 41, 20, 10, 5, 2) \qquad \quad l_{331} = 7   \quad \qquad \frac{t_{331}+2l_{331} - 1}{331} < \frac{1+7}{2^7}< 0.1
\]

\begin{align*}
&p=104,729 \qquad \quad  LH(104,729) = (52364, 26182, 13091, 6545, 3272, 1636, 818, 409, 204, 102, 51, 15, 12, 6, 3)\\
 & \qquad  \qquad l_{104,729}=15  \qquad \qquad \frac{t_{104,729}+2l_{104,729} - 1}{104,729} < \frac{1+15}{2^{15}}< 0.0005
\end{align*}

\bigbreak

\section{Preview of future work.} \label{sect:preview}

\noindent

In this Section we offer the reader a preview of our next project, applying the Teneva Game to rational knots, with the following graphical calculations of minimum numbers of colors. We refer the reader to \cite{kl-k} and \cite{kldmtcs} for further information and terminology on rational knots.

\subsection{The rational knot R(5, 2)}

\noindent

\begin{figure}[!ht]
	\psfrag{0}{\huge$0$}
	\psfrag{1}{\huge$1$}
	\psfrag{2}{\huge$2$}
	\psfrag{3}{\huge$3$}
	\psfrag{4}{\huge$4$}
	\psfrag{5}{\huge$5$}
	\psfrag{6}{\huge$6$}
	\psfrag{7}{\huge$7$}
	\psfrag{8}{\huge$8$}
	\psfrag{9}{\huge$9$}
	\psfrag{10}{\huge$10$}
	\psfrag{11}{\huge$11$}
	\psfrag{12}{\huge$12$}
	\psfrag{13}{\huge$\bf 13$}
	\psfrag{8-12 removed}{\huge$8-12 \text{ removed}$}
	\psfrag{k1-1 colors}{\huge$\text{i.e., } k_1-1 \text{ colors removed}$}
	\psfrag{5 removed}{\huge$5 \text{ removed}$}
	\psfrag{k2-2 colors}{\huge$\text{i.e., } k_2-2 \text{ colors removed}$}
	\psfrag{5 removed}{\huge$5 \text{ removed}$}
	\psfrag{1 colors}{\huge$\text{i.e., } 1 \text{ color removed}$}
	\psfrag{k1=6}{\huge$k_1=6$}
	\psfrag{k2=3}{\huge$k_2=3$}
	\psfrag{b2}{\huge$b_2$}
	\psfrag{b3}{\huge$b_3$}
	\psfrag{b4}{\huge$b_4$}
	\psfrag{b5}{\huge$b_5$}
	\psfrag{b6}{\huge$b_6$}
	\centerline{\scalebox{.5}{\includegraphics{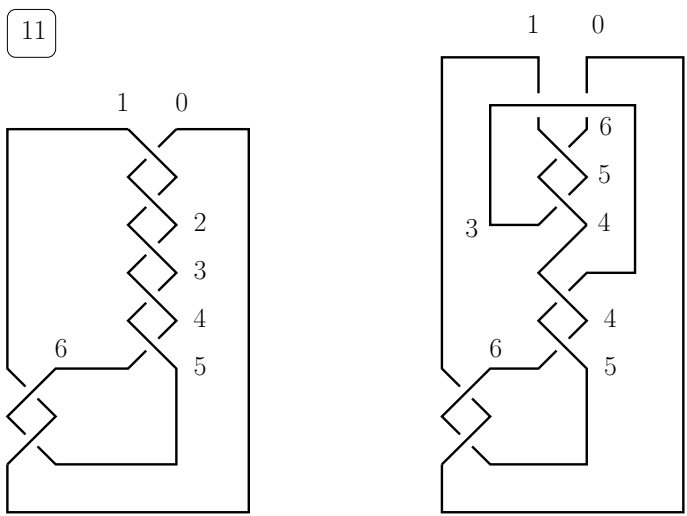}}}
	\caption{$R(5, 2)$ admits a non-trivial $11$-coloring with $6$ colors.}\label{fig:r52}
\end{figure}

The leftmost diagram of Figure \ref{fig:r52} depicts a rational knot, $R(5, 2)$.
The fraction associated to this rational knot is
\[
\frac{1}{5+\frac{1}{2}} = \frac{2}{11}
\]
so its determinant is $11$.

\begin{figure}[!ht]
	\psfrag{0}{\huge$0$}
	\psfrag{1}{\huge$1$}
	\psfrag{2}{\huge$2$}
	\psfrag{3}{\huge$3$}
	\psfrag{4}{\huge$4$}
	\psfrag{5}{\huge$5$}
	\psfrag{6}{\huge$6$}
	\psfrag{7}{\huge$7$}
	\psfrag{8}{\huge$8$}
	\psfrag{9}{\huge$9$}
	\psfrag{10}{\huge$10$}
	\psfrag{11}{\huge$11$}
	\psfrag{12}{\huge$12$}
	\psfrag{13}{\huge$\bf 13$}
	\psfrag{8-12 removed}{\huge$8-12 \text{ removed}$}
	\psfrag{k1-1 colors}{\huge$\text{i.e., } k_1-1 \text{ colors removed}$}
	\psfrag{5 removed}{\huge$5 \text{ removed}$}
	\psfrag{k2-2 colors}{\huge$\text{i.e., } k_2-2 \text{ colors removed}$}
	\psfrag{5 removed}{\huge$5 \text{ removed}$}
	\psfrag{1 colors}{\huge$\text{i.e., } 1 \text{ color removed}$}
	\psfrag{k1=6}{\huge$k_1=6$}
	\psfrag{k2=3}{\huge$k_2=3$}
	\psfrag{b2}{\huge$b_2$}
	\psfrag{b3}{\huge$b_3$}
	\psfrag{b4}{\huge$b_4$}
	\psfrag{b5}{\huge$b_5$}
	\psfrag{b6}{\huge$b_6$}
	\centerline{\scalebox{.5}{\includegraphics{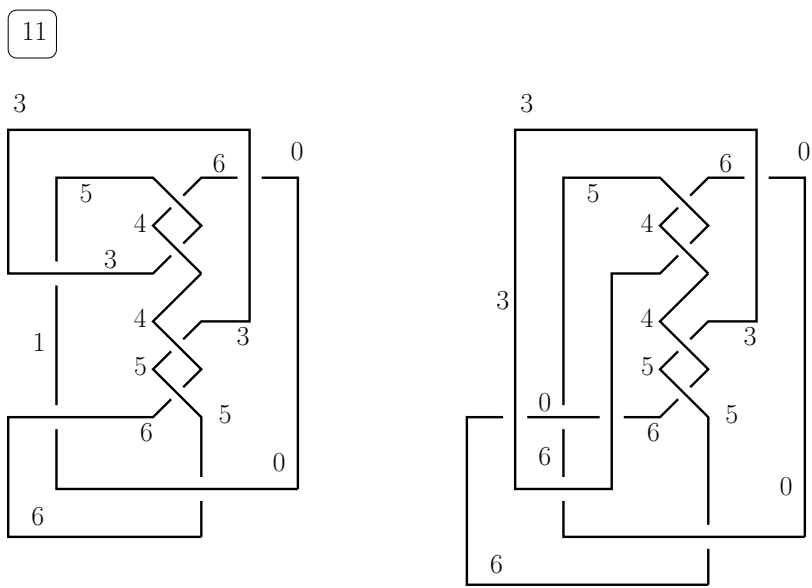}}}
	\caption{$R(5, 2)$ admits a non-trivial $11$-coloring with $5$ colors: rearranging the rightmost diagram of Figure \ref{fig:r52} and a move to eliminate color $1$.}\label{fig:r52cont}
\end{figure}

Inspection of Figures \ref{fig:r52} and \ref{fig:r52cont} show that $R(5, 2)$ admits a non-trivial $11$-coloring with $5$ colors. Since the determinant of $R(5, 2)$ is $11$, then thanks to the aforementioned result in \cite{lm},
\[
mincol_{11} R(5, 2) = 5.
\]

\bigbreak

\subsection{The rational knot R(4, -3)}

\noindent

\begin{figure}[!ht]
	\psfrag{0}{\huge$0$}
	\psfrag{1}{\huge$1$}
	\psfrag{2}{\huge$2$}
	\psfrag{3}{\huge$3$}
	\psfrag{4}{\huge$4$}
	\psfrag{5}{\huge$5$}
	\psfrag{6}{\huge$6$}
	\psfrag{7}{\huge$7$}
	\psfrag{8}{\huge$8$}
	\psfrag{9}{\huge$9$}
	\psfrag{10}{\huge$10$}
	\psfrag{11}{\huge$11$}
	\psfrag{12}{\huge$12$}
	\psfrag{13}{\huge$\bf 13$}
	\psfrag{8-12 removed}{\huge$8-12 \text{ removed}$}
	\psfrag{k1-1 colors}{\huge$\text{i.e., } k_1-1 \text{ colors removed}$}
	\psfrag{5 removed}{\huge$5 \text{ removed}$}
	\psfrag{k2-2 colors}{\huge$\text{i.e., } k_2-2 \text{ colors removed}$}
	\psfrag{5 removed}{\huge$5 \text{ removed}$}
	\psfrag{1 colors}{\huge$\text{i.e., } 1 \text{ color removed}$}
	\psfrag{k1=6}{\huge$k_1=6$}
	\psfrag{k2=3}{\huge$k_2=3$}
	\psfrag{b2}{\huge$b_2$}
	\psfrag{b3}{\huge$b_3$}
	\psfrag{b4}{\huge$b_4$}
	\psfrag{b5}{\huge$b_5$}
	\psfrag{b6}{\huge$b_6$}
	\centerline{\scalebox{.5}{\includegraphics{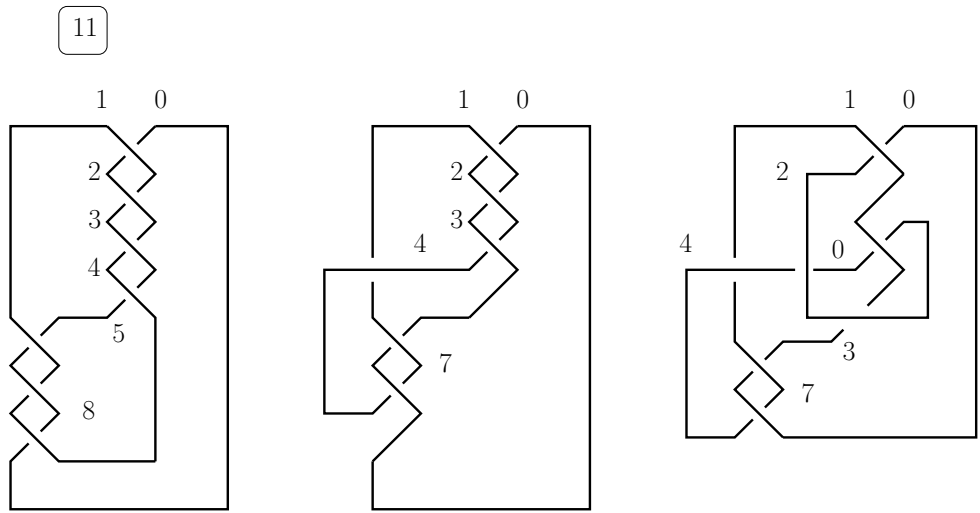}}}
	\caption{$R(4, -3)$ admits a non-trivial $11$-coloring with $6$ colors.}\label{fig:r43m}
\end{figure}

The leftmost diagram of Figure \ref{fig:r43m} depicts a rational knot, $R(4, -3)$.
The fraction associated to this knot is
\[
\frac{1}{4+\frac{1}{-3}} = \frac{-3}{11}
\]
so its determinant is $11$.

\begin{figure}[!ht]
	\psfrag{0}{\huge$0$}
	\psfrag{1}{\huge$1$}
	\psfrag{2}{\huge$2$}
	\psfrag{3}{\huge$3$}
	\psfrag{4}{\huge$4$}
	\psfrag{5}{\huge$5$}
	\psfrag{6}{\huge$6$}
	\psfrag{7}{\huge$7$}
	\psfrag{8}{\huge$8$}
	\psfrag{9}{\huge$9$}
	\psfrag{10}{\huge$10$}
	\psfrag{11}{\huge$11$}
	\psfrag{12}{\huge$12$}
	\psfrag{13}{\huge$\bf 13$}
	\psfrag{8-12 removed}{\huge$8-12 \text{ removed}$}
	\psfrag{k1-1 colors}{\huge$\text{i.e., } k_1-1 \text{ colors removed}$}
	\psfrag{5 removed}{\huge$5 \text{ removed}$}
	\psfrag{k2-2 colors}{\huge$\text{i.e., } k_2-2 \text{ colors removed}$}
	\psfrag{5 removed}{\huge$5 \text{ removed}$}
	\psfrag{1 colors}{\huge$\text{i.e., } 1 \text{ color removed}$}
	\psfrag{k1=6}{\huge$k_1=6$}
	\psfrag{k2=3}{\huge$k_2=3$}
	\psfrag{b2}{\huge$b_2$}
	\psfrag{b3}{\huge$b_3$}
	\psfrag{b4}{\huge$b_4$}
	\psfrag{b5}{\huge$b_5$}
	\psfrag{b6}{\huge$b_6$}
	\centerline{\scalebox{.5}{\includegraphics{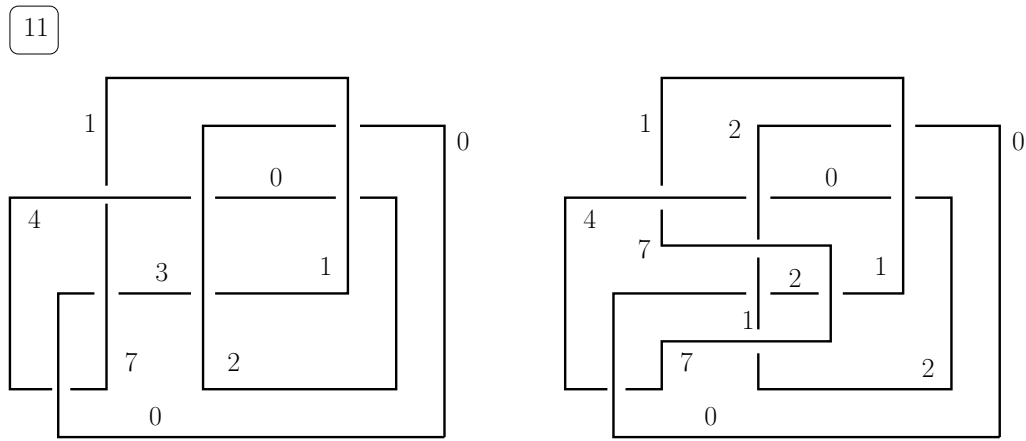}}}
	\caption{$R(4, -3)$ admits a non-trivial $11$-coloring with $5$ colors: rearranging the rightmost diagram of Figure \ref{fig:r43m} and a move to eliminate color $3$.}\label{fig:r43mcont}
\end{figure}

Inspection of Figures \ref{fig:r43m} and \ref{fig:r43mcont} show that $R(4, -3)$ admits a non-trivial $11$-coloring with $5$ colors.
Arguing as with the preceding rational knot
\[
mincol_{11}R(4, -3) = 5.
\]

\bigbreak
\section{Acknowledgements.} \label{sect:ackn}

\noindent

P.L. acknowledges support from FCT (Funda\c c\~ao para a Ci\^encia e a Tecnologia), Portugal, through project number PTDC/MAT/101503/2008, ``New Geometry and Topology''. P.L. also thanks the School of Mathematical Sciences at the University of Nottingham for hospitality during his stay there where the present article was written.

L. K. and P.L. thank Shin Satoh for Corollary \ref{cor:mincol}.


\end{document}